\documentclass[leqno,12pt]{amsart} 
\usepackage{amssymb} 
\usepackage[english]{babel}
\usepackage{amsthm}
\usepackage{amsmath}
\usepackage{amssymb}
\usepackage[latin1]{inputenc}
\def\R{\text{$\mathbb{R}$}}

\def\pr{\parallel}

\def\lra{\longrightarrow}

\def\lra{\longrightarrow}

\def\d1#1#2{\frac{d#1}{d#2}}

\def\p1#1#2{\frac{\partial #1}{\partial #2}}
\def\to{t_o}
\def\pr{\parallel}

\def\part{a=\to \leq t_1 \leq \ldots \leq t_{n-1} \leq t_{n}=b}

\def\P{\text{$\mathbb{P}$}}
\def\C{\text{$\mathbb{C}$}}
\def\CP{\text{$\mathbb{CP}$}}
\def\P{\text{$\mathbb{P}$}}
\newcommand\SO{{\rm SO}}

\newcommand\SU{{\rm SU}}

\def\Sp{{\rm Sp}}
\newcommand\U{{\rm U}}

\newcommand\mf[1]{\mathfrak{#1}}

\newcommand\To{{\rm T}}
\def\o{\omega}

\def\mr{ M_{\lambda}}

\theoremstyle{plain} 
\newtheorem{thm}{\indent\sc Theorem}[section]

\newtheorem{cor}[thm]{\indent\sc Corollary}
\newtheorem{prop}[thm]{\indent\sc Proposition}

\theoremstyle{definition} 

\newtheorem{rem}[thm]{\indent\sc Remark}
\newtheorem{exe}[thm]{\indent\sc Example}

\begin{document}
\title[Homogeneous Lagrangian submanifolds]{Hamiltonian actions 
and homogeneous Lagrangian submanifolds} 

\author[L. Biliotti]{Leonardo Biliotti} 
\dedicatory{Dedicated to Professor Francesco Mercuri on his sixtieth birthday}

%
\address{%
Dipartimento di Matematica \endgraf
Universit\`a Politecnica delle Marche \endgraf
Via Brecce Bianche\endgraf
60131, Ancona Italy
}
\email{biliotti@dipmat.univpm.it}

%
\begin{abstract}
We consider a connected symplectic manifold $M$ 
acted on properly and in a Hamiltonian fashion 
by a connected Lie group $G$. 
Inspired to the recent paper \cite{gb2}, see also \cite{ch} and \cite{pacini}, 
we study Lagrangian orbits of Hamiltonian actions. The dimension of
the moduli space of the Lagrangian orbits is given and we also 
describe under which condition a Lagrangian orbit is isolated. 
If $M$ is a compact K\"ahler manifold we give a necessary and sufficient 
condition  to an isometric action admits a Lagrangian orbit.
Then we investigate Lagrangian homogeneous submanifolds
on the symplectic cut and on the  symplectic reduction. 
As an application of our results, we  give new examples 
of Lagrangian homogeneous submanifolds 
on the blow-up at one point
of the complex projective space and on the weighted projective spaces.
Finally, applying Proposition \ref{slice} that
we may call \emph{Lagrangian slice theorem} for group acting with
a fixed point, we give new examples of Lagrangian homogeneous
submanifolds on irreducible Hermitian symmetric spaces of compact and 
noncompact type.
\end{abstract}
\maketitle
\section{Introduction} 
Let $(M,\omega)$ be a symplectic manifold. A Lagrangian submanifold of $M$
is a submanifold of half dimension of $M$ on which the symplectic form
$\omega$ vanishes. Lagrangian submanifolds are intensively studied and  they
have classically played an important role in symplectic geometry
(see \cite{hl}, \cite{ta}, \cite{oh} \cite{oh2}). Recently their role
has been expanded beyond that of their use in understanding 
symplectic diffeomorphisms.

In \cite{oh} the author asks for a group theoretical machinery
producing Lagrangian subma\-ni\-folds in Hermitian symmetric spaces
and in a recent paper \cite{gb2} the existence problem of
homogeneous Lagrangian submanifolds in compact K\"ahler manifolds
is studied, obtaining a characterization of isometric actions
admitting a Lagrangian orbit for a large class of compact K\"ahler
manifolds including irreducible Hermitian symmetric spaces.

In this paper we first study Lagrangian homogeneous submanifolds in a
symplectic ma\-ni\-fold. We give the dimension of the moduli space of 
Lagrangian orbits,
describing under which condition a Lagrangian orbit of 
a reductive Lie group $G$ is isolated.   
The uniqueness result generalizes Theorem $2$ in \cite{gb2}.

Our main tool will be the moment map that can be 
defined whenever we consider Hamiltonian action on $M$. More precisely,
let $(M,\omega)$ be a connected symplectic manifold,
$G$ be a connected Lie group of  symplectic diffeomorphisms  acting
in a Hamiltonian fashion. This means  there exists a 
map $\mu:M \lra \mf g^*$, where $\mf g$ is the Lie algebra
of $G$, which is called \emph{moment map}, satisfying:
\begin{enumerate}
\item for each $X\in \mf{g}$ let
\begin{itemize}
\item $\mu^{X}: M \lra \R,\ \mu^{X} (p)= \mu(p) (X),$ the
component of $\mu$ along $X,$ \item $X^\#$ be the vector field on
$M$ generated by the one para\-me\-ter subgroup $\{ \exp (tX):
t \in \R \} \subseteq G$.
\end{itemize}
Then
$$
{\rm d} \mu^{X}= {\rm i}_{X^\#} \o ,
$$
i.e. $\mu^{X}$ is a Hamiltonian function for the vector field
$X^\#.$
\item $\mu$ is $G-$equivariant, i.e. $\mu (gp)=Ad^* (g)
(\mu(p)),$ where $Ad^*$ is the coadjoint representation on $\mf g^*$.
\end{enumerate}
In general the matter of existence/uniqueness of $\mu$ is
delicate. However whenever $\mf g$ is semisimple the moment map
exists and it is unique \cite{gs}. If $(M,\omega)$ is a compact
K\"ahler manifold and $G$ is a connected compact Lie group of
holomorphic isometries then the existence problem is resolved (see
\cite{HW}): a moment map exists if and only if $G$ acts trivially
on the Albanese torus Alb($M$). 

In the sequel we always assume that for
every $\alpha \in \mf g^*$, $G\alpha$ 
is a locally closed coadjoint orbit of $G$. 
Observe that the condition of a coadjoint orbit being locally
closed is automatic for reductive group and for their semidirect
products with
vector spaces. There exists an example of a solvable Lie group due to
Mautner \cite{mau} p.$512$, with non-locally closed coadjoint
orbits. These assumptions are needed to applying the symplectic
slice, see \cite{lb}, \cite{gs}, \cite{or}, \cite{sl}, 
and the symplectic stratification of the reduced space given
in \cite{lb}.

Before we give the statement of our first main result, we fix our notation.

Let $Gx=G/G_x$ be a $G$-orbit. Since $G_x$ is compact, we may split
$\mf g=\mf g_x \oplus \mathbf m$, as $G_x$-modules.
We denote by  $\mf{n (g_{{\it x}})}$ the Lie algebra of $N(G_x)$ the normalizer
of $G_x$ in $G$. Let $v=v_x + v_m \in \mf g_x \oplus
\mathbf m$ be an element of $\mf{n( g_{{\it x}} )}$.  
Then $[v_m, \mf g_x] \subset \mf g_x$, 
i.e. $v_m \in \mf{n(g_{{\it x}} )}$, 
which implies that $[v_m , \mf g_x]=0$,
since $\mathbf m$ is $G_x$-invariant. This means
\begin{equation} \label{poi}
\mf{n (g_{{\it x}})} = \mf g_x \oplus \{ v\in \mathbf m: [v,\mf g_x]=0
\}=\mf g_x \oplus \mathbf s .
\end{equation}
Let $\mf{z(g)}$ be the Lie algebra of the center of $G$. Clearly
$\mf{z(g)}\subseteq \mf{n(g_{{\it x}})}$ and the projection 
$\pi:\mf{n(g_{{\it x}})}\lra \mf g_x$ maps $\mf{z(g)}$ to 
$\mf{z(g_{{\it x}}) }$.
\begin{thm} \label{lagr}
Let $(M,\o )$ be a connected symplectic $G$-Hamiltonian manifold
with moment map $\mu: M\lra \mf g$. Assume that $G/G_x=Gx$ is a
Lagrangian orbit. Then the dimension of the moduli space of the
Lagrangian orbits which contains $Gx$ 
is $\dim \mf {z(g)} \cap \mathbf s$. Therefore
the dimension of the moduli space of the Lagrangian
orbits is equal or less than $\dim N(G_x)/G_x$ and $Gx$ is an
isolated Lagrangian orbit, i.e. there exists a $G$-invariant
neighborhood on which if $Gy$ is Lagrangian then $Gy=Gx$, if and
only if the projection $\pi:\mf{n(g_{{\it x}})}\lra \mf g_x$
maps $\mf{z(g)}$ one-to-one to $\mf{z(g_{{\it x}})}$. 
Hence if $G$ is a semisimple Lie group, a Lagrangian
orbit, if there exists, is isolated. 
Moreover, at any level
set $\mu^{-1}(c)$  if there exists a Lagrangian orbit, it is
isolated.
\end{thm}
Note that we have proved under which condition an action have infinitely
many Lagrangian orbits. 
The next result characterizes ones having isolated 
Lagrangian orbit, generalizing Theorem $2$ in \cite{gb2}.
\begin{thm} \label{ls}
Let $G$ be a reductive Lie group acting properly and in a Hamiltonian
fashion on a connected symplectic manifold $M$. A Lagrangian $G$-orbit, $Gp$
is isolated if and only if there exists a semisimple closed
Lie subgroup $G'$ of $G$
such that $Gp=G'p$.
\end{thm}
We may also characterize isometric actions on 
compact K\"ahler manifold admitting a Lagrangian orbit, applying
a result of Kirwan \cite{ki} and the symplectic stratification of the 
reduced space given in \cite{lb}.
\begin{thm} \label{zio}
Let $G$ be a compact connected Lie group acting isometrically and 
in a Hamiltonian fashion on a compact K\"ahler manifold $M$. 
If we denote by $\mu$ the correspon\-di\-ng 
moment map, then $G$ admits a Lagrangian orbit if and only if there exists 
$c \in \mf{z(g)}$ such that $\mu^{-1}(c)$ is a Lagrangian submanifold.  
\end{thm}
As an immediately corollary we have the following result.
\begin{cor}
If $G$ is a compact Lie group acting isometrically
in a Hamiltonian fashion on compact 
K\"ahler manifold then at any level set $\mu^{-1}(c)$ there exists at most one
Lagrangian orbit. Moreover, if $G$ is a compact semisimple Lie group, 
then $G$ admits a Lagrangian orbit if and only if $\mu^{-1}(0)$ is a Lagrangian 
submanifold of $M$.
\end{cor}
If $G$ is compact, it is standard fix an $Ad(G)$-invariant scalar product 
$\langle \cdot,\cdot \rangle$ and we identify $\mf g$ with $\mf g^*$ by
means $\langle \cdot,\cdot \rangle$, we can think $\mu$ as a 
$\mf g$-valued map. It is also natural study the squared moment map 
$\parallel \mu \parallel^2$ and its critical set. This function has been 
intensively studied in \cite{ki}  obtaining strong information on the 
topology of $M$. In \cite{Bi3} we have proved that if a point $x$ realizes a 
local maximum of $\parallel \mu \parallel^2$ then the $G$-orbit through $x$ is
symplectic. It is natural to study the ``dual''  problem: 
the points which realize the minimum of $\parallel \mu \parallel^2$. Note that 
a Lagrangian orbit could describe this set by Theorem \ref{zio} whenever $M$ 
is a compact K\"ahler manifold.

Next, we use Hamiltonian actions to construct non-standard 
homogeneous Lagrangian submanifolds on K\"aheler manifolds. 
The first ``non-standard'' examples on 
the complex projective space appeared in \cite{ch}. More recently,
see \cite{gb2}, the classification of isometric actions of simple Lie groups
admitting a Lagrangian orbit on complex projective spaces is given.     

We study Lagrangian homogeneous submanifolds on the
symplectic reduction (\cite{cannas}) and on the symplectic cut
(see \cite{gbl} and \cite{lerman})
making a connection between Lagrangian orbits on $M$ and Lagrangian orbits
on symplectic reduction and on symplectic cut, see Proposition \ref{su}
and Proposition \ref{cut}.
These results are interesting since the complex 
projective space is the reduced space of the standard ${\rm S}^1$-action
on $\C^n$ and the symplectic cut can be obtained by blowing-up $M$ along a 
symplectic submanifold,   
see \cite{gbl} and \cite{lerman} for more details.

We use these results to construct non-standard  
homogeneous Lagrangian submanifolds on the blow-up at one point of the
complex projective space (example \ref{exlag} and subsection \ref{scoppia}).
and we also deduce, using the classification given in \cite{gb2}, see
Corollary \ref{class} and Remark \ref{lc},
the classification of the simple compact Lie groups $K$ acting 
isometrically on  $\C^n$ such that
${\rm S}^1 \times K$ admits a Lagrangian orbit on $\C^n$.

Since  our result holds   also when the symplectic reduction is an orbifold,
see \cite{cannas} for more detail about orbifold, in subsection \ref{pw}
we give new examples of homogeneous Lagrangian submanifolds on 
weighted projective spaces .

Finally,  applying Proposition \ref{slice} 
which we may call \emph{Lagrangian slice}
for group ac\-ti\-ng with a fixed point, we give new examples of Lagrangian
homogeneous submanifolds  on
irreducible Hermitian symmetric spaces of compact and noncompact type.
\begin{prop} \label{lagherm}
Let $G$ be a group appears in Table $1$. Then $G$ admits a Lagrangian orbit.
\end{prop}
\begin{small}
\[
\begin{array}{|c|c|c|}\hline
G & M  & M^*\\ \hline\hline
{\rm G}_2 \subseteq \SO(7) \subseteq \SO(8) & \SO(8) / \U(4)  & \\ \hline
{\rm SO}(2) \times \SO(n),\ n\geq3 & \SO(n+2)/ \SO(2) \times \SO(n) 
& {\rm SO}_o (2,n) / \SO(2) \times \SO(n) \\ \hline
{\rm S}(\U(1) \times \U(n)) & \CP^n & 
\SU(1,n)/ {\rm S}(\U(1) \times \U(n))\\ \hline
Z({\rm S}(\U(2) \times \U(2n))) \times \Sp(n),\ n\geq2&
\SU(2n+2)/ {\rm S}(\U(2) \times \U(2n)) 
& \SU (2,2n)/ {\rm S}(\U(2) \times \U(2n))\\ \hline
\U(2n)   &  \SU(4n) / \U(2n) & \SU^*(4n)/\U(2n) \\ \hline
\U(2n+1)\subset \U(2n+2),\ n\geq 2 & \SU(4n+4))/ \U(2n+2)) & 
\SU^* (4n+4)/\U(2n+2)\\ \hline
\U(n) & \Sp(n)/ \U(n) & \Sp(n,\R)/ \U(n)\\ \hline
\To^1 \cdot {\rm E}_6 &  {\rm E}_7/ \To^1 \cdot {\rm E}_6  
& {\rm E}^{-25}_7/ \To^1 \cdot {\rm E}_6\\ \hline
\To^1 \cdot {\rm Spin}(9) \subseteq \To^1 \cdot {\rm Spin}(10)&
{\rm E}_6 / \To^1 \cdot {\rm Spin}(10) & 
{\rm E}^{-14}_6 / \To^1 \cdot {\rm Spin}(10)\\ \hline
{\rm Spin}(7)\subset {\rm SO}(8) \subset \SU(8) & 
\SU(8) / {\rm S}(\U(2) \times \U(6)) & \\
\hline
\end{array}
\]
\end{small}
\section{Existence and uniqueness}
Here we follows the notation as in \cite{lb} and in the introduction.
The first easy remark is the following: if $Gx$ is a
Lagrangian submani\-fold of $M$ then $\mu(x)\in \mf{z(g)}^*$. 
Indeed, since
$Ker {\rm d}\mu_x=
(T_x G(x))^{\perp_{\omega}}$
we have $\mu_{|_{Gx}}=c\in \mf g^*$ which
implies that $c\in \mf{z(g)}^*$.\\
$\ $ \\
\emph{Proof of Theorem \ref{lagr}}. Suppose $Gx$ is a Lagrangian
orbit. We may assume that $\mu(x)=0$.
From the symplectic slice, there exists a $G$-invariant
neighborhood which is symplectomorphic to $Y= G \times_{G_x} (\mf g /
\mf g_x)^*$ and the moment map is given by
\[
\mu([g,v])=Ad(g)^*(j(v)). 
\]
We recall that we may split
$\mf g=\mf g_x \oplus \mathbf m$  as $G_x$-modules and
$j$ is induced by the above decomposition,
see \cite{lb}. It is well-known, shrinking the neighborhood of
the zero section of $Y$ if necessary, that $\dim G[g,v] \geq \dim
Gx$. Hence $G[e,v]$ is Lagrangian if and only if
$\mu([e,v]) \in \mf{z(g)}^*$ if and only if $v\in \mathbf m^* \cap
\mf{z(g)}^*$ if and only if  $v\in \mf{z(g)}^* \cap \mathbf s^*$
by (\ref{poi}). This claims that the dimension of the 
moduli space of Lagrangian orbits which
contains $Gx$ is $\dim (\mf{z(g)}\cap \mathbf s) \leq \dim
N(G_x)/G_x$, from (\ref{poi}). 
We also deduce that $Gx$ is isolated if and only if 
$\mf{z(g)}\cap \mathbf s =\{0 \}$ if and only if the projection
of $\mf g$ into $\mf g_x$ maps $\mf{z(g)}$ one-to-one to 
$\mf{z(g_{{\it x}})}$.  
Finally,  since in a $G$-invariant neighborhood of a Lagrangian orbit, 
the moment map is given by $\mu([g,m])=Ad^*(g)(j(m))$, the moment map locally
separates $G$-orbits, concluding our proof.
\begin{flushright}
$\square$
\end{flushright}
\begin{exe} \label{qui}
We consider $G=\To^1 \times \SU(2)$
acting on $\C^2 \oplus \C^2$
as follows
\[
(t,A)(v,w)=(At^{-1}v,Aw)
\]
This action is Hamiltonian with moment map
\begin{small}
\[
\mu((z_1,z_2),(w_1,w_2))= \frac{i}{2}
\left( \begin{array}{cc}
\parallel w_1 \parallel^2 -
\parallel w_2 \parallel^2 & z_1\overline{z_2}+ w_1\overline{w_2} \\
z_2 \overline{z_1}+ w_2\overline{w_1} 
& \parallel w_2 \parallel^2 -\parallel w_1 \parallel^2 
\end{array}
\right) + \frac{i}{2} \parallel ( z_1 , z_2 )  \parallel^2 .
\]
\end{small}
The orbit through $((1,1),(1,-1))$ is Lagrangian  and the $\To^1$-orbit through
$((1,1),(1-1))$ induces  a curve in the moduli space of the Lagrangian orbits.
\end {exe}
\emph{Proof of Theorem \ref{ls}.} 
In the sequel we denote by $G_{{\rm}ss}$ the closed Lie semisimple subgroup
of $G$ whose Lie algebra is $[\mf g, \mf g]$. 

Assume there exists a closed
semisimple Lie subgroup $G'$ of $G$ 
such that $Gp=G'p$ is Lagrangian. 
It is well-known, see \cite{Path}, there exists a $G$-invariant neighborhood
$Y$ of $Gp$ such that $\dim Gq \geq \dim Gp$ for every $q\in Y$. 
Since $G'\subset G_{{\rm ss}}$ we get $G_{{\rm ss}}p=Gp$. If we denote
by $\overline\mu$ the moment map corresponding to the 
$G_{{\rm ss}}$-action on $M$, one may note that $\overline\mu(q)=0$
if and only if $\mu(q)\in \mf{z(g)}$, where $\mu$ is the moment map for the 
$G$-action on $M$, since 
the moment map of the
$G_{ss}$-action is the projection on $\mf g_{ss}$ of $\mu$.
Therefore, shrinking $Y$ if necessary, given $q\in Y$ 
we get that $Gq$ is Lagrangian if and only if $G_{{\rm ss}}q$ is, 
proving $Gp$ is isolated, since $G_{{\rm ss}}$ is semisimple.

Vice-versa, let $Gp$ be an isolated Lagrangian orbit. 
We may assume $\mu(p)=0$. Since $Gp$ is isolated, 
there exists a $G$-invariant neighborhood $Y$ such that
$\mu(Y) \cap \mf{z (g)}=\{ 0 \}$. Hence, given  $x \in Y$ we have
$\overline{\mu}(x)=0$ if and only if $\mu(x)\in \mf{z(g)}$ if and
only if $\mu(x)=0$. This proves $Gp=\mu^{-1}(0) \cap Y= \overline
{\mu}^{-1}(0) \cap Y$.

Let $H$ be a principal isotropy to the $G_{{\rm ss}}$-action on
$Gp$. It is well-known, see \cite{lb} and \cite{sl}, that the
intersection of the stratum $M^{(H)}$ of orbit type $(H)$ with the
zero level set of the moment map $\overline{\mu}$ is a submanifold
of $M$ of constant rank and the orbit space
\[
(M_o )^{(H)}= ( \overline{\mu}^{-1}(0) \cap M^{(H)} )/G_{ss},
\]
has a natural symplectic structure $(\omega_o )_{(H)}$ whose
pullback to $\overline{\mu}^{-1}(0) \cap M^{(H)}$ coincides with
the restriction to $\overline{\mu}^{-1}(0) \cap M^{(H)}$ of the
symplectic form on $M$. Since $\overline{\mu}^{-1}(0) \cap Y$ is
Lagrangian, we get $(\omega_o )_{(H)}=0$ which claims
$(M_o)^{(H)}$ are points. This implies $G_{{\rm ss}}$ acts
transitively on $Gp$ concluding our proof.
\begin{flushright}
$\square$
\end{flushright}
\emph{Proof of the Proposition \ref{zio}}.
Assume that $G$ admits a Lagrangian orbit. We may suppose that $Gp$ is 
Lagrangian and $\mu(p)=0$, where $\mu$ is the corresponding moment map. 
It is easy to check that $G^{\C}p$ is open in $M$. Therefore $M$ is projective
algebraic and there exists a $G$-equivariant embedding in some projective space
\cite{HW}. In particular, since the $G^{\C}$-action on $M$ is algebraic, 
$G^{\C}p$ is Zariski open in $M$.

If $x\in \mu^{-1}(0)$ lies in different $G$-orbit, then  there exist, due 
a results of Kirwan \cite{ki}, two $G^{\C}$-invariant disjoint neighborhood
$U_p$ and $U_x$ of $p$ and $x$ respectively. 
Since $U_p$ contains $Gp$, the neighborhood 
$U_p$ must meet $U_x $,  which is an absurd.
We point out that, this claim that there exists
at most one Lagrangian orbit at any level set $\mu^{-1}(c)$.

Vice-versa, assume that $\mu^{-1}(c)$ is Lagrangian. 
Let $H$ be a isotropy subgroup for the $G$-action on $M$.
Since 
$M^{(H)} \cap \mu^{-1}(c)/G$ has a natural symplectic structure $(\o )_H$ whose
pullback to $M^{(H)} \cap \mu^{-1}(c)$ coincides with the restriction to
$M^{(H)} \cap \mu^{-1}(c)$ of the symplectic form on $M$, we get that 
$M^{(H)} \cap \mu^{-1}(c)/G$ are points. 
This claim that given $p\in \mu^{-1}(c)$ we have 
$\mu^{-1}(c)=Gp$, since $\mu^{-1}(c)$ is connected, concluding our proof.  
\begin{flushright}
$\square$
\end{flushright}
\begin{rem}
In the sequel we always assume that $G$ acts properly and 
coisotropically on $M$. This means that there exists
an open dense subset $U$ such that for every $x\in U$, $Gx$ is a
coisotropic submanifold of $M$ with respect
to $\o$, i.e.  $(T_x Gx )^{\perp_{\o}} \subset T_x Gx$.
Coisotropic actions are intensively studied in \cite{HW} and \cite{GS2}
and in the following paper
\cite{PT},\cite{bg}, \cite{Bi} the complete classification
of compact connected Lie groups acting isometrically and
in a Hamiltonian fashion on irreducible compact Hermitian symmetric
spaces are given. In a forthcoming paper \cite{Bi2}, we shall study
coisotropic actions  of Lie groups acting properly and in a 
Hamiltonian fashion on a symplectic manifold $M$.
An equivalent condition for a connected Lie group $G$ acts  
coisotropically on $M$ is that
for every $\alpha \in \mf g^*$ the set $G\mu^{-1}(\alpha)/G$ are points 
\cite{Bi2}. In
particular if the fibers of the moment map $\mu:M \lra \mf g^*$
are connected, this holds if $M$ and $G$ are  compact, then $G$ admits a
Lagrangian orbit if and only if $\mu^{-1}(z)$ is a Lagrangian submanifold
for some $z\in \mf{z(g)}^*$.

Now, assume that a principal orbit is Lagrangian.
Then there exists an abelian closed Lie group $T$ which have the same
orbits of the $G$-action on $M$. Indeed, let $H$ be a principal isotropy. 
Since $\mu (M^{(H)}) \subset \mf{z(g)}$, $M$ is mapped by $\mu$ to 
$\mf {z(g)}$. Therefore, given $x\in M^{H}$,
from  symplectic slice, we have that the
$H$-submodule $\mathbf m$ such that $\mf g=\mf h \oplus \mathbf m$
is abelian. Let $T$ be the closure of the torus whose Lie algebra
is $\mathbf m$. Note that $Tx=Gx$. Therefore, the set of regular points
of $T$ is a subset of the set of $G$-regular points which implies
that $T$ and $G$ have the same orbits, since both the $T$-action and the 
$G$-action are proper actions.
\end{rem}
\section{Symplectic cut, symplectic reduction
and homogeneous Lagrangian submanifolds} 
In this section we shall
investigate how the existence of homogeneous Lagrangian
submanifolds on $M$ implies the existence of homogeneous
Lagrangian submanifold on the symplectic cut and on the symplectic
reduction. For sake of completeness we describe briefly these
constructions and we invite the reader to see \cite{cannas},
\cite{lerman} and \cite{gbl} for a good exposition of these
subjects.

Let $M$ be  a connected symplectic manifold and let $G$ be a
compact connected Lie group acting in a Hamiltonian fashion on $M$. 
Let $K$ be a semisimple compact Lie subgroup of $G$ and let $\To^k$ be a
$k$-dimensional torus which centralizes $K$ in $G$, i.e. $\To^k
\subset C_G (K)$. In the sequel we denote by
\[
\phi:M \lra \mf k \oplus  \mf t_k,
\]
where $\mf t_k=$Lie$(\To^k)$, the moment map of the $\To^k \cdot K$-action
on $M$ and with $\mu$, respectively with $\psi$, 
the moment map of the $K$-action on $M$, respectively 
a moment map of the $\To^k$-action on $M$.

Let $\lambda \in \mf t_k$ be such that
$\To^k$ acts freely on $\psi^{-1}(\lambda)$. Then the symplectic reduction
\[
M_{\lambda}= \psi^{-1}(\lambda) / \To^k ,
\]
is a symplectic manifold on which $K$ acts, since it commutes with
$\To^k$, in a Hamiltonian fashion with moment map
\[
\overline{\mu}: M_{\lambda} \lra \mf k , \ \overline\mu ([x])= \mu(x).
\]
This proves that $\overline\mu([x])= 0$ if and only if
$\mu(x)=0$ if and only if $\phi(x)\in \mf t_k$.

Let $[p]\in M_{\lambda}$. It is easy to see that $k[p]=[p]$
if and only if there exists $r(k)\in \To^k$,
which is unique  since $\To^k$ acts freely on $\psi^{-1}(\lambda)$,
such that $kp=r(k)p$.
This means that the following homomorphism of Lie groups
\[
K_{[p]} \stackrel{F}{\lra} (\To^k \cdot K)_p,\ F(k)=k r(k)^{-1},
\]
is a covering map, due the fact that $K$ is semisimple. Hence
\[
\dim K[p]= \dim (\To^k \cdot K)p - \dim \To^k.
\]
Since $\dim M= \dim M_{\lambda} - 2\dim \To^k$, we
have proved the following result.
\begin{prop} \label{su}
Let $M$ be a connected symplectic manifold on which a compact connected 
Lie group $G$ acts in a Hamiltonian fashion on $M$. 
Let $K$ be a compact semisimple Lie subgroup of $G$ and let
$\To^k$ be a $k$-dimensional torus which centralizes $K$ in $G$. Let 
$\lambda \in \mf t_k$ be such  that 
$\To^k$ acts freely on $\psi^{-1}(\lambda)$, where $\psi$ is a moment map of 
the $\To^k$-action on $M$. Then $(\To^k \cdot K)p$, $p\in \psi^{-1}(\lambda)$, 
is Lagrangian if and only if $K[p]$ in $M_{\lambda}$ is.
\end {prop}
\begin{rem} \label{pp}
We would like to point out that Proposition \ref{su} holds
if we assume only that $\lambda \in$Lie($\To^k$) is a regular value.
In this case the symplectic reduction could be an orbifold, see \cite{cannas},
and $\To^k$ could act almost freely on the level set $\psi^{-1}(\lambda)$.
Hence the following
map
\[
F : K_{[p]} \lra (\To^k \cdot K)_p / \To_p^k ,
\]
is a covering map which implies, since $\To_p^k$ is finite,
$\dim K_{[p]}=\dim (\To^k \cdot K)_p$.
\end{rem}
It is well-known that the complex projective space is the
K\"ahler reduction for the standard ${\rm S}^1$-action on 
$\C^{n+1}$, see \cite{cannas}.
Therefore, from Proposition \ref{su}, we deduce the following result.
\begin{cor} \label{class}
Let $K$ be a compact semisimple Lie group acting
in a Hamiltonian fashion  on $\CP^n$.
Then $K$ admits a Lagrangian orbit if and only if the
${\rm S}^1\times K$-action on $\C^{n+1}$ admits.
\end{cor}
\begin{rem} \label{lc}
From the classification given in \cite{gb2},
we deduce which compact simple Lie groups $K$ satisfy 
$K \times{\rm S}^1$ act with a Lagrangian orbit on $\C^n$. 
\end{rem}
Now, we consider a one-dimensional torus $\To^1$ which
centralizes a connected compact Lie group $K$. We may consider the
symplectic cut which is briefly described in the following (see
also \cite{lerman} and \cite{gbl}). In the sequel we think the moment maps
for the $\To^1$-actions as $\R$-valued maps.  

We consider the symplectic manifold $M \times \C$ with symplectic form
$\omega  -\frac{i}{2}dz \wedge d\overline{z}$.
$\To^1$ acts diagonally on $M \times \C$ as $t(m,z)=(tm,tz)$ with moment map
$\overline{\psi}= \psi + \parallel z \parallel^2$.

Let $\lambda \in \R=$Lie$(\To^1)$ be such that $\To^1$
acts freely on $\psi^{-1}(\lambda)$. Then $\To^1$ acts freely on
$\overline{\psi}^{-1}(\lambda)$ so we may consider the symplectic reduction,
\[
M^{\lambda}= \overline{\psi}^{-1}(\lambda)/ \To^1 ,
\]
which is called symplectic cut. Note that $\dim M= \dim M^{\lambda}$.
$K$ acts on $M^{\lambda}$ as
$k[m,z]=[km,z]$ with moment map
\[
\overline{\mu}([x,z])=\mu(x) , 
\]
where $\mu$ is the moment map of the $K$-action on $M$. Since
if $z \neq 0$, $K_{[m,z]}=K_{m}$ we deduce  that $Km$ is
Lagrangian if and only if $K[m,z]$ is.
\begin{prop} \label{cut}
Let $(M,\omega)$ be a symplectic manifold on which a compact connected
Lie group $G$ acts in a Hamiltonian fashion. Assume that $K$ is a compact
subgroup of $G$  whose centralizer in $G$ contains a one-dimensional torus 
${\rm T}^1$.
Then $K$ admits a Lagrangian orbit in the open subset 
$\{x\in M: \psi(x)<\lambda\}$, where 
$\psi$ is the corresponding moment map of the $\To^1$-action on $M$, 
if and only if 
the $K$-action on the symplectic cut, 
obtained from the $\To^1$-action on $M$, admits.
\end{prop}
\begin{exe} \label{exlag}
Let $\To^1$ acting on $\CP^n$ as follows
$$
(t, [z_o, \ldots, z_n])  \longrightarrow [t^{-1}z_o,z_1, \ldots ,z_n] .
$$
This is a Hamiltonian action with moment map
$$
\phi ([z_o, \ldots, z_n] )= \frac{1}{2} \frac{\pr z_o \pr^2}{\pr
z_o \pr^2 + \cdots + \pr z_n \pr^2} \ .
$$
One can prove that $\phi([1, \ldots , 0])$ is the global maximum
of the moment map $\phi$ and $\phi^{-1} (\frac{1}{2})=[1, \ldots ,
0].$ Hence, as in \cite{gbl} page 5, if $\lambda=\frac{1}{2}-
\epsilon,$ $\epsilon \cong 0,$ then the K\"ahler cut $(\CP^n
)^{\lambda}$ is the blow-up of $\CP^n$ at 
$[1, \ldots , 0]$ that we shall indicate with $\CP^n_{[1,\ldots,0]}$.
The torus action $\To^n$ on $\CP^n$ given by
\[
(t_1, \ldots, t_n)([z_o,\ldots,z_n]=[z_o, t^{-1}_1 z_1, \ldots , t^{-1}_n z_n],
\]
is Hamiltonian and the principal orbits are Lagrangian.
Note that $\To^n$ acts in a Hamiltonian fashion on $\CP_{[1, \ldots,0]}^n$ and
for the above discussion we conclude that the $\To^n$-action on
$\CP_{[1, \ldots,0]}$ admits  a Lagrangian orbit. Moreover,
it acts coisotropically and its principal orbits are Lagrangian.

Let $\To^1=Z({\rm S}(\U(1) \times \U(4))$ and
$K=\To^1 \times \U(2)$ acting on $\CP^4$ as follows
\[
(t,A)(Z)=(t,A)[z_o,z_1,z_2,w_1,w_2]=(t,A)[z_o,z,w]=
[tz_o, A t^{-4}z, At^{-4}w].
\]
This action is Hamiltonian with moment map
\begin{small}
\[
\mu(Z)=
\frac{i}{2 \parallel Z \parallel^2 }
\left( \begin{array}{cc}
\parallel z \parallel^2  + \parallel w
\parallel^2 & z_1\overline{z_2}+ w_1\overline{w_2} \\
z_2 \overline{z_1}+ w_2\overline{w_1} &
\parallel z \parallel^2  + \parallel w \parallel^2
\end{array}
\right) +
\frac{i}{2 \parallel Z \parallel^2}
( -\parallel z_o \parallel^2 + \frac{1}{2}
(\parallel (z,w) \parallel^2))
\]
\end{small}
Let $p=[z_o,1,1,1,-1] \in \CP^4$. Note that $\mu(p)
\in$Lie($\To^1)\oplus \mf{z(u(2))}$, and if $z_o \neq 0$ then
$\dim Kp=4$, which implies $Kp$ is Lagrangian.

Since $K$ commutes with the above $\To^1$-action, we obtain that the
$K$-action on $\CP_{[1,\ldots,0]}^4$ admits a Lagrangian orbit, which is 
$K[p]$.
\end{exe}
Next, we claim the Lagrangian slice theorem for $
G$-action with a fixed point.
\begin{prop} \label{slice}
Let $G$ be a compact Lie group acting in a Hamiltonian fashion with
a fixed point on a symplectic manifold $M$. 
If the slice representation at the fixed point
has a Lagrangian orbit then so has the $G$-action on $M$.
\end{prop}
\begin{proof}
It follows immediately from symplectic slice. Indeed, if $Gp=p$ then
$\mu(p)=\beta \in \mf{z(g)}$ and from symplectic slice, 
the moment map is locally given by
\[
\mu(gp,m)=\beta +Ad^*(g)(\mu_V (m))
\]
If $Gm$ is a Lagrangian orbit of the $G$-action on the slice then
$\dim Gm=\frac{1}{2} \dim M$ and $\mu_V(m)\in \mf{z(g)}$.
In particular
$\dim G[p,m]=\frac{1}{2} \dim M$ and $\mu([p,m]) \in \mf{z(g)}$
which implies $G[p,m]$ is Lagrangian.
\end{proof}
\emph{Proof of the Proposition \ref{lagherm}}.
Except for the first and the last cases, 
the prove follows from the classification given in
\cite{gb2}, Corollary \ref{class} and finally
from Proposition \ref{slice}.  In the sequel we always refer to 
Table $1$ p. 16 in \cite{gb2}, which has been enclosed at the end of this
section, 
for compact simple Lie groups acting with a Lagrangian orbit in $\CP^n$.
We briefly explain our method. 

The semisimple part of $G$ admits a
Lagrangian orbit on the projective space of the slice, by Table $1$. 
Hence, by Corollary \ref{su}, $G$ has a Lagrangian orbit on the slice which
implies, from Proposition \ref{slice}, $G$ admits a Lagrangian
orbit on $M$.  
We consider only the Hermitian symmetric spaces of compact type
since for the non compact case we have only to change $M$ 
with $M^*$.
\begin{enumerate}
\item ${\rm G}_2$ acting on $M=\SO(8)/\U(4)$.
We use the same argument as in \cite{Bi}.  Since 
${\rm G}_2 \cap \U(4)= \SU(3)$,
the orbit through $[\U(4)]$ is Lagrangian. Indeed, let 
\[
\phi:M \lra \mf g^*_2
\]
be the moment map of the ${\rm G}_2$-action on $M$. 
One may check that \[
\dim {\rm G}_2[\U(4)]=\frac{1}{2}\dim \SO(8)/\U(4),
\] 
${\rm G}_2 \phi([\U(4)])={\rm G}_2 /{\rm P}$ is a flag manifold and
$\SU(3) \subseteq {\rm P}$. Since $\SU(3)$ is a maximal subgroup of
${\rm G}_2$ which does not centralize a torus, we deduce that $P={\rm G}_2$
proving $\phi([\U(4)])=0$. Hence ${\rm G}_2 [\U(4)]$ is Lagrangian.
\item $G=\SO(2) \times \SO(n)$ acting on $M=\SO(n+2) / \SO(2) \times \SO(n)$.  
The slice is $\C^n$ on which $\SO(n)$ acts with $\Lambda_1$.
Since $\SO(n)$ admits a Lagrangian orbit on $\CP^{n-1}$,
from Corollary \ref{class}, 
$G$ admits a Lagrangian orbit on the slice which implies, 
from Proposition \ref{slice}, $G$ admits on $M$;
\item $G={\rm S}(\U(1) \times \U(n))$ acting on $M=\CP^n$. Since
$\SU(n)$ has a Lagrangian orbits on $\CP^{n-1}$ then $G$ admits a
Lagrangian orbit on the slice which implies that $G$ has a
Lagrangian orbit on $M$.
\item
$G=Z({\rm S}(\U(2) \times \U(2n))) \times \Sp(n)$ acting on
$M=\SU(2n+2)/ {\rm S}(\U(2) \times \U(2n))$. The slice is given by
$\C^{2n} \oplus \C^{2n}$ on which $\Sp(n)$ acts diagonally while the one
dimensional torus acts as
\[
t(v,w)=(t^{\frac{1-n}{n} } v, t^{\frac{1-n}{n}} w)
\] 
Since $\Sp(n)$ admits a Lagrangian orbit on $\CP^{4n-1}$, we get 
that $G$ admits on $M$.
\item
$G=\U(2n)$ acting on  $M=  \SU(4n) / \U(2n)$. The slice is given by
 $\Lambda^2 (\C^{2n})$ and $\SU(2n)$ acts with a Lagrangian orbit on
its complex projective space. Hence $G$ admits a Lagrangian orbit
on $M$.
\item $G=\U(2n+1)$ acting on  $M=\SU(4n+4)/ \U(2n+2)$. The slice is given by
\[
\Lambda^2 (\C^{2n+1}) \oplus \C \otimes \C^{2n+1}.
\]
Note that the center acts as
\[
t(X,v)=(t^2X,tv).
\]
As in \cite{gb2}, one may prove that $Gp$ is Lagrangian, where
\[
p=( \left(\begin{array}{cc}
0& 0\\
0 & J_n
\end{array}\right)
,e_1 ),
\] 
and
\[
J_n=\left(\begin{array}{cc}
0   & -I_n \\
I_n &  0 \\
\end{array}\right),
\]
which implies that $G$ admits a Lagrangian orbit on $M$;
\item
$G=\U(n)$ acting on $ \Sp(n)/ \U(n)$. The slice is given by
${\rm S}^2(\C^n)$ and $\SU(n)$ has a Lagrangian orbit on the
projective space of
${\rm S}^2(\C^n)$. Then $G$ admits a Lagrangian orbit on $M$.
\item 
$G=\To^1 \cdot {\rm E}_6$ acting on
$M={\rm E}_7/ \To^1 \cdot {\rm E}_6$. The slice is given by
$\C^{27}$ on which $G$ acts with the $\Lambda_1$ representation.
Since ${\rm E}_6$ admits a Lagrangian orbit on $\CP^{26}$, $G$ admits on $M$.
\item $G=\To^1 \cdot {\rm Spin}(9) \subseteq \To^1 \cdot {\rm Spin}(10)$
acting on $ {\rm E}_6 / \To^1 \cdot {\rm Spin}(10)$. The slice is given by
$\C^{16}$ on which $G$ acts with the \emph{spin-representation}.
Since ${\rm Spin}(9)$ has a Lagrangian orbit on $\CP^{15}$, $G$
admits a Lagrangian orbit on $M$ from Corollary \ref{class} and 
Proposition \ref{slice}.
\item This case  is proved in \cite{bg} p. $1736$.
\end{enumerate}
\begin{flushright}
$\square$
\end{flushright}
\subsection{Homogeneous Lagrangian submanifolds on weighted
projective Spaces} \label{pw}
In \cite{gb2} it was proved that the following action 
\[
\begin{array}{|c|c|c|c|}\hline
G & \rho  & \dim_{\C} \P(V)& G^o_p \\ \hline\hline
\SU(n) & \Lambda_1 \oplus \Lambda^*_1 & 2n-1& \SO(n) \\ \hline
\SU(2n+1)& \Lambda_2\oplus \Lambda_1 & 2n^2+3n+1& \Sp(n) \\ \hline
\Sp(n) & \Lambda_1 \oplus \Lambda_1 & 4n-1 & \Sp(n-1) \\ \hline
{\rm Spin}(10) & \Lambda_e \oplus \Lambda_e & 31 & \SU(5) \\ \hline
\end{array}
\]
admit Lagrangian orbit. We briefly explain the notation, which follows one
given in \cite{gb2}. $\rho$ denotes the representation and $G^o_p$ denotes
the connected component of the identity of the 
isotropy of the Lagrangian orbit. Moreover we have identified
the fundamental weights $\Lambda_i$ with the corresponding irreducible 
respresentations. 

We shall prove that if $G$ appears in the above table,
then it induces an action on some weighted projective space with a 
Lagrangian orbit. 

Let $\To^1$ be a one-dimensional torus acting on $V=V_1 \oplus V_2$ as
\[
t(v,w)=(t^{-k}v,t^{-s}),
\]  
where $k,s$ are distinct natural numbers. This actions is Hamiltonian with
moment map
\[
\psi((v,w))=\frac{i}{2}(k \parallel v \parallel^2 + s \parallel w \parallel^2).
\]
Let $\lambda=\frac{i}{2}$. The reduced space $\mr$ is the weighted
projective space, which we denote by $\P(V)_{[k,s]}$, on which $G$ acts
in a Hamiltonian fashion since it commutes with the $\To^1$-action. 
The map
$\mu:V \lra \mf{g}^*$ defined for
every $(v,w)\in V$ and for every $X\in \mf{g}$ by
\[
\mu((v,w))=-i \langle X(v,w),(v,w)\rangle
\] 
is the moment map for the $G$-action on $V$, where
$\langle \cdot , \cdot \rangle$ denotes the natural hermitian scalar product
on $V$. 

We claim that $L=G \times \To^1$ admits a Lagrangian
orbit which lies in $\psi^{-1}(\frac{i}{2})$. Note that the corresponding
moment map for the $L$-action is $\xi=\mu + \psi$. Our approach is the following. 

Let
$p=(v,w)\in V$ be such that $G[p]$ is Lagrangian in $\P(V)$, which has been
calculated in \cite{gb2}. Then we may prove that $L\overline{p}$, where
$\overline{p}= \frac{1}{
\sqrt{k \parallel v \parallel^2 + s \parallel w \parallel^2}}p$,
is Lagrangian, which implies that $G[\overline p ]$ is Lagrangian in 
$\P(V)_{[k,s]}$.
\begin{enumerate}
\item $G=\SU(n)$ and $\rho=\Lambda_1 \oplus \Lambda^*_1 $.  
$p=(e_1,e^*_1)$ and one may prove that the $L$-orbit through $\overline p$
has dimension $2n$. Since $\xi(\overline p )=\psi(\overline p )$, 
$L \overline p$ is 
Lagrangian;
\item $G=\SU(2n+1)$ and $\rho=\Lambda_2 \oplus \Lambda_1$. 
$p=(J_n,e_1)$, where $J_n$ is the same as in case $(6)$
of the proof of Proposition \ref{slice}. 
Since $\dim L \overline p =2n^2+3n+2$, $L \overline p $ is Lagrangian;
\item $\Sp(n)$ and $\rho=\Lambda_1\oplus \Lambda_1$. 
$p=(e_1,e_2 )$ and it is easy to check that 
$\dim L \overline p =4n$, which implies that $L \overline p$ is Lagrangian;
\item $G={\rm Spin}(10)$ and $\rho=\Lambda_e \oplus \Lambda_e$.
$p=(1+e_{1234},e_{15}+e_{2345})$ (see \cite{sk} for the notation). 
One may prove that  $\dim L \overline p=32$ which means that $L$ admits a
Lagrangian orbit. 
\end{enumerate}    
\subsection{Homogeneous Lagrangian submanifolds on the blow-up at 
one point of the complex projective space} \label{scoppia}
Let $K$ be compact Lie group acting linearly on $V$ with at two 
sub-modules,  
i.e. $V=V_1\oplus V_2$ and $K$ preserves $V_i$, $i=1,2$. Suppose that
$K$ admits a Lagrangian orbit on the complex projective space
$P(V)$. If we consider a one-dimensional torus $\To^1$ acting on
$V_1$ or on $V_2$ which commutes with $K$ 
we may induces an Hamiltonian action of $K$ on the K\"ahler cut 
$\P(V)^{\lambda}$ and this action admits
a Lagrangian orbit. In \cite{gb2} it was proved that the following action 
\[
\begin{array}{|c|c|c|c|}\hline
G & \rho  & \dim_{\C} \P(V)& G^o_p \\ \hline\hline
\SU(n) & \Lambda_1 \oplus \Lambda^*_1 & 2n-1& \SO(n) \\ \hline
\SU(2n+1)& \Lambda_2\oplus \Lambda_1 & 2n^2+3n+1& \Sp(n) \\ \hline
\Sp(n) & \Lambda_1 \oplus \Lambda_1 & 4n-1 & \Sp(n-1) \\ \hline
{\rm Spin}(10) & \Lambda_e \oplus \Lambda_e & 31 & \SU(5) \\ \hline
\end{array}
\]
admits a Lagrangian orbit.

We shall prove that these group admits a Lagrangian
orbit on the blow-up at one point of $\P (V)$. We analyze in detail only
the first case: the other ones are similar. 

Let $\To^1$ be a torus acting on $\P(\C^n \oplus (\C^n )^*)$ 
as $t[(v,w)])=[( t^{-1}v,w)]$. 
This action is Hamiltonian with moment map
\[
\mu([(v,w)])=\frac{i \parallel v \parallel^2}{2\parallel [v,w] \parallel^2} .
\]
Hence $p=[(1,\ldots,1),(0\ldots,0)]$ is the global  
maximum and as usual, see \cite{gbl},
given $\lambda=\frac{1}{2} -\epsilon$ $\epsilon \cong 0$, the K\"ahler-cut
is the blow-up of the complex projective space $\P(V)$ at $p$, 
which we denote by
$\P(\C^n \oplus (\C^n )^* )_{[p]}$. Since  the $\SU(n)$-action on
$\Lambda_1\oplus \Lambda^*_1$ commutes with the $\To^1$-action, from
Proposition \ref{cut}, $\SU(n)$ admits a Lagrangian orbit on 
$\P(\C^n \oplus (\C^n )^* )_{[p]}$.    
Summing up we have the following Homogeneous Lagrangian submanifolds.
\[
\begin{array}{|c|c|c|}\hline
G & \P(V)_{[p]}& G^o_p \\ \hline\hline
\SU(n) & \P(\C^n \oplus (\C^n )^*)_{[(1,\ldots,1),(0,\ldots ,0)]} 
& \SO(n) \\ \hline
\SU(2n+1)& \P(\Lambda^2(\C^{2n+1})\oplus \C^{2n+1} )_{[(0,\ldots,0),(1,\ldots,1)]} 
& \Sp(n) \\ \hline
\Sp(n) & \P(\C^{2n} \oplus \C^{2n})_{[(1,\ldots,1),(0,\ldots,0)]} 
& \Sp(n-1) \\ \hline
{\rm Spin}(10) 
& \P(\C^{16} \oplus \C^{16} )_{[(1,\ldots,1),(0,\ldots,0)]} 
& \SU(5) \\ \hline
\end{array}
\]
\newpage
\begin{center}
\textbf{Table:} Lagranagian orbits of simple Lie groups in projective spaces 
\end{center}
$\ $ \\
\[
\begin{array}{|c|c|c|c|}\hline
G & \rho & \dim_{\C} \P(V)&{\rm cond.} \\ \hline\hline
\SU(n) & 2\Lambda_1 & \frac{n(n+1)}{2}-1 & \\ \hline
\SU(n) & \Lambda_1 \oplus \Lambda^*_1 & 2n-1 & \\ \hline
\SU(n) & \underbrace{\Lambda_1 \oplus \cdots \oplus \Lambda_1}_{n} 
& n^2 -1 & \\ \hline
\SU(2n) & \Lambda_2 & n(2n-1)-1 & n\geq 3 \\ \hline
\SU(2n+1) & \Lambda_2 \oplus \Lambda_1 & 2n^2 +3n +1 & n\geq 2 \\ \hline
\SU(2) & 3 \Lambda_1 & 3 & \\ \hline
\SU(6) & \Lambda_3 & 19 & \\ \hline
\SU(7) & \Lambda_3 & 34 & \\ \hline
\SU(8) & \Lambda_3 & 55 & \\ \hline
\Sp(n) & \Lambda_1 \oplus \Lambda_1 & 4n-1 & \\  \hline
\Sp(3) & \Lambda_3 & 13 & \\ \hline
\SO(n) & \Lambda_1 & n-1 & n\geq 3 \\ \hline
{\rm Spin}(7) & {\rm spin. rep.} & 7 & \\ \hline
{\rm Spin}(9) & {\rm spin. rep.} & 15 & \\ \hline
{\rm Spin}(10) & \Lambda_e \oplus \Lambda_e & 31 & \\ \hline
{\rm Spin}(11) & {\rm spin. rep.} & 31 & \\ \hline
{\rm Spin}(12) & \Lambda_e & 31 & \\ \hline 
{\rm Spin}(14) & \Lambda_e & 63 & \\ \hline
{\rm E}_6 & \Lambda_1 & 26 & \\ \hline
{\rm E}_7& \Lambda_1 & 55 & \\ \hline
{\rm G}_2 & \Lambda_2 & 6 & \\ \hline
\end{array}
\]
\newpage


\begin{thebibliography}{99}
\bibitem{alq}
\textsc{H. Azad, J. Loeb and M. Qureshi},
\newblock \emph{Totally real orbits in affine quotients of reductive groups,}
\newblock{Nagoya Math. J.} \textbf{139}, (1983) 87--92.
%
\bibitem{lb}
\textsc{L. Bates and E. Lermann},
\newblock \emph{Proper group actions and symplectic stratified spaces },
\newblock{Pacific J. Math.} \textbf{181}, (1997) 201--229.
%
\bibitem{gb2}
\textsc{L. Bedulli and Gori},
\newblock \emph{Homogeneous Lagrangian submanifold}
\newblock arXiv:math.DG/0604169;
%
\bibitem{besse}
\textsc{A. L. Besse},
\newblock \emph{Einstein manifold}, Springer-Verlang (1978).
%
\bibitem{Bi3}
{\sc L. Biliotti},
\newblock \emph{A note on the moment map on symplectic manifolds,}
\newblock{preprint arXiv math./DG\-06\-05\-2\-70}
%
\bibitem{bg}
\textsc{L. Biliotti and A. Gori},
\newblock \emph{Coisotropic and polar actions on complex {G}rassmannians},
\newblock {Trans. Amer. Math. Soc.} \textbf{357} (2005), 
no. \textbf{5}, 1731--1751.
\bibitem{Bi}
\textsc{L. Biliotti}
\newblock \emph{Coisotropic and polar actions on compact irreducible
 {H}ermitian symmetric spaces},
\newblock Trans. Amer. Math. Soc. \textbf{358} (2006) 3003--3022.
%
\bibitem{Bi2}
\textsc{L. Biliotti},
\newblock \emph{Some results on multiplicity-free spaces}, preprint.
%
\bibitem{gbl}
\textsc{D. Burns, V. Guillemin and E. Lerman},
\newblock \emph{K\"ahler cuts},
\newblock arXiv: math.DG/\-0212062.
%
\bibitem{cannas}
\textsc{A. da Silva Cannas},
\newblock{\emph Lectures on Symplectic Geometry,}
\newblock{Lecture Notes in Math. \textbf{1764},
Springer-Verlag, Berlin, 2001.}
%
%
\bibitem{CE}
{\sc Cheeger, J., and Ebin, D.},
\newblock \emph{Comparison theorems in {R}iemannian geometry},
\newblock North--Holland Publishing Co., Amsterdam, 1975.
%
\bibitem{ch}
\textsc{ R. Chiang},
\newblock \emph{New Lagrangian submanifolds of} $\CP^n$,
\newblock Int. Math. Res. Not. \textbf{45}, (2004) 2437-2441.
%
%
%
%
%
\bibitem{gs}
\textsc{V. Guillemin and S. Stenberg},
\newblock \emph{Symplectic techniques in physics}, Cambridge Univ. Press,
Cambridge, 1990.
%
\bibitem{GS2}
\textsc{V. Guillemin and S. Stenberg},
\newblock \emph{Multiplicity-free spaces},
\newblock{J. Differential Geom. {\textbf 19} (1984), 31--56}.
%
\bibitem{hl}
\textsc{R. Harvey and H. B.Lawson},
\newblock \emph{Calibrates Geometries},
\newblock{Acta Math.} \textbf{148}, (1982) 47--157.
%
%
\bibitem{he}
\textsc{S. Helgason},
\newblock \emph{Differential Geometry, Lie groups, and Symmetric spaces}
\newblock Academic Press New-York-London, second edition (1978).
%
\bibitem{HW}
\textsc{ A. Huckleberry and T. Wurzbacher},
\newblock \emph{Multiplicity-free complex manifolds},
\newblock Math. Ann. \textbf{286}, (1990) 261--280.
%
%
%
\bibitem{ki}
\textsc{Kirwan},
\newblock \emph{Cohomology of quotients in symplectic and algebraic geometry},
\newblock Math. Notes \textbf{31} Princeton (1984).
%
%
\bibitem{KN}
\textsc{S. Kobayashi and N. Nomizu},
\newblock \emph{Foundations of Differential
Geometry,} \textbf{Vol. I}, Interscience Publisher, J. Wiles \&
Sons, New-York, 1963.
%
\bibitem{lerman}
\textsc{E. Lerman}
\newblock\emph{Symplectic cut}
\newblock{Math. Res. Lett. \textbf{2}, (1995) 247--258}
%
\bibitem{oh}
\textsc{Y.-G. Oh},
\newblock \emph{Mean Curvature vector of symplectic topology of Lagrangian
submanifolds in K\"ahler-Einstein manifolds,}
\newblock{Math. Z.} \textbf{216}, (1990) 471--482.
%
\bibitem{oh2}
\textsc{Y.-G. Oh},
\newblock \emph{Second variation and stabilities of minimal
Lagrangian subma\-ni\-folds in K\"ahler ma\-ni\-folds,}
\newblock{Invent. Math.} \textbf{101}, (1990) 501--519.
%
\bibitem{or}
\textsc{J.P. Ortega and T. S. Ratiu},
\newblock \emph{A symplectic slice Theorem,}
\newblock{Lett. Math. Phys. {\textbf 58} 81--93 (2002).}
%
\bibitem{pacini}
\textsc{T. Pacini},
\newblock \emph{Mean Curvature flow, orbits, moment maps},
\newblock Trans. Amer. Math. Soc. \textbf{355} (2004) 3343--3357.
%
\bibitem{Path}
\textsc{R. S. Palais and C. L. Terng},
\newblock \emph{Critical point theory and submanifold
  geometry},
\newblock Lecture Notes in Math. \textbf{1353}, Springer-Verlag, New-York, 
1988.
%
\bibitem{PT}
\textsc{F. Podest\`a and  G. Thorbergsson},
\newblock \emph{Polar and coisotropic actions on {K\"ahler} manifolds},
\newblock {Trans.Amer. Math. Soc.} \textbf{354} (2002), 1759--1781.
%
\bibitem{sk}
\textsc{M. Sato and T. Kimura,}
\newblock \emph{A classification of irreducible prehomogeneous vector spaces
and their relative invariants,}
\newblock Nagoya Math.  J. \textbf{65}, 1--155 (1977).
%
\bibitem{sl}
\textsc{R. Sjamaar,  and E. Lerman},
\emph{Stratified symplectic spaces and reduction},
Ann. of Math. {\textbf 134}, (1991) 375--422.
%
%
\bibitem{mau}
\textsc{L. Pukanszky,}
\newblock \emph{Unitary representations of solvable groups},
\newblock {Ann. Sci. \'Ecole Normale Sup.} \textbf{4} (1972), 475--608.
%
\bibitem{ta}
\textsc{M. Takeuchi},
\newblock \emph{Stability of certain minimal submanifolds of
compact Hermitian symmetric spaces,}
\newblock{Tohoku J. Math.} \textbf{(2)} \textbf{36}, (1984) 293--314.
\end{thebibliography}
\end{document}